\documentclass [11pt]{article}
 \usepackage{amssymb}
 \usepackage{amsmath,amsthm}
 \usepackage{mathrsfs}
 \usepackage{color}
\usepackage[numbers,sort&compress]{natbib}
 \usepackage{graphicx}
\usepackage{tikz}
 \usepackage[colorlinks=true]{hyperref}

\normalsize

\newcommand\Bin{\operatorname{Bin}}
\renewcommand{\dfrac}[2]{\lower0.15ex\hbox{\large$\textstyle\frac{#1}{#2}$}}

\newtheorem{theorem}{Theorem}[section]
\newtheorem{lemma}[theorem]{Lemma}

\newtheorem{remark}[theorem]{Remark}

\numberwithin{equation}{section}

\begin{document}
\title
{\bf Random $K_k$-removal algorithm 
\thanks{The work was partially supported by NSFC.}}
\author
{{\large Fang Tian$^1$\thanks{Corresponding Author:\ tianf@mail.shufe.edu.cn(Email Address).}\quad Zi-Long Liu$^2$\quad Xiang-Feng Pan$^3$}\\
{\small $^1$ Department of Applied Mathematics}\\%
{\small Shanghai University of Finance and Economics, Shanghai, 200433, China} \\
{\small\tt tianf@mail.shufe.edu.cn}\\[1ex]
{\small $^2$School of Optical-Electrical and Computer Engineering}\\
{\small University of Shanghai for Science and Technology, Shanghai,
200093, China}\\
{\small\tt liuzl@usst.edu.cn}\\[1ex]
{\small $^3$School of Mathematical Sciences}\\
{\small  Anhui University, Hefei, Anhui, 230601,  China}\\
{\small\tt xfpan@ahu.edu.cn}}
\date{}
 \maketitle

\begin{abstract}
One interesting question is how a graph develops from some constrained random graph process,
which  is a fundamental mechanism in the formation and evolution of dynamic networks.
The problem here is referred to
the random $K_k$-removal algorithm.
For a fixed integer $k\geqslant 3$, it starts with
a complete graph on $n\rightarrow\infty$ vertices and iteratively removes the edges of an
uniformly chosen $K_k$. This algorithm terminates once no $K_k$s remain
and at the same time it generates one linear $k$-uniform hypergraph.
For $k=3$, it was shown that the size in the final graph is $n^{3/2+o(1)}$.
Less  results are on the cases when $k\geqslant 4$.
In this paper, we prove that the exact expected trajectories of various key parameters in the algorithm
to some iteration such that the final size in the algorithm
is at most  $n^{2-1/(k(k-1)-2)+o(1)}$ for $k\geqslant 4$.
We also show the bound is a natural barrier.
\end{abstract}

\hskip 10pt{\bf Keywords:}\ random greedy algorithm, $K_k$-free,
the critical interval method,  dynamic concentration.

\hskip 10pt {\bf Mathematics Subject Classifications:}\ 05D40, 68R10

\section{Introduction}

Extremal problems are central research issues in random graph algorithms,
 which are also fundamental mechanisms in the formation and evolution of dynamic networks.
 A better understanding of the underlying graph offers us opportunities to study
how a graph develops from some constrained random greedy process.
Recently, the power of random
greedy algorithm is illustrated in~\cite{guo2021} by showing the existence
of mathematical objects with better properties. Each time random
greedy algorithms go beyond classical applications of the probabilistic method
used in previous work.

 The problem here is referred to 
 the random $K_k$-removal algorithm.
 Given a fixed integer $k\geqslant 3$, the random $K_k$-removal algorithm for generating
one $K_k$-free graph, and at the same time creating a linear $k$-uniform hypergraph,
is defined as follows. Start from a complete graph on vertex set $[n]$,
denoted by $G(0)$, and $G(i+1)$ is the remaining graph from
$G(i)$ by selecting one $K_k$ uniformly at random out of all $K_k$s in $G(i)$
and deleting all its edges. Let the hitting time $M$ be $M=\min\{i\,:G(i)\ {\mbox{is\ } K_k{\mbox {-free}}}\}$
and $E(i)$ denote the edge set of $G(i)$, thus $|E(M)|$ is the number of edges
in the final $K_k$-free graph.

Work on finding the exact values of $|E(M)|$ has
evolved over the past 20 years and is a nontrivial task even for $k=3$.
Bollob\'{a}s and Erd\H{o}s~\cite{bollobas98}
conjectured that with high probability
$|E(M)|=n^{3/2+o(1)}$ when $k=3$. It was shown $|E(M)|=o(n^2)$
by Spencer~\cite{spencer95} and independently by R\"{o}dl and Thoma~\cite{rt96}.
Grable~\cite{grable97} improved this bound to $|E(M)|\leqslant n^{7/4+o(1)}$.
Bohman et al.~\cite{bohman101} introduced the critical interval method
for proving dynamic concentrations.  They~\cite{bohman15}
confirmed the exponent in a breakthrough by
 generalizing the approach in~\cite{bohman101}.
Less results directly studied the random $K_k$-removal algorithms when $k\geqslant 4$.
Bennett and Bohman~\cite{bennett15} conjectured that $|E(M)|\leqslant n^{2k/(k+1)+o(1)}$
as a folklore for $k\geqslant 3$ when they
investigated the  random greedy hypergraph matching
algorithm. It is exactly the one proposed by Bollob\'{a}s and Erd\H{o}s when $k=3$.

 A different recipe for obtaining a random $K_k$-free graph is the so-called ``$K_k$-free process".
 In that algorithm, starting with an empty graph, the ${n\choose 2}$ edges are randomly inserted
so long as no $K_k$s are formed in the current graph.
Despite the high similarity between the two protocols,
 it was shown~\cite{bohman15} that the random $K_k$-removal
algorithm has proved quite challenging at the level of acquiring
the correct exponent of the final number of edges.
A pseudo-random heuristic for divining the evolution of various key
parameters plays a central role in the understanding of these algorithms
that produce interesting combinatorial objects 
~\cite{bennett15,bohman09,bohman15,bohman101,bohman19,pico142,war142}.

In this paper,  we directly discuss the structure of
random $K_k$-removal algorithm for $k\geqslant 4$.
We design an ensemble of appropriate random variables
including the number of $K_k$s,
using a heuristic assumption to find the trajectories of
these variables when the process evolves.
Compared with the random $K_3$-removal algorithm, it is  challenging
to make use of these  auxiliary variables to  analyze the one-step
change of the number of $K_k$s when $k\geqslant 4$ and
show a rigorous proof of their expressions. At last, we verified that
 \begin{theorem}
Given a fixed integer $k\geqslant 4$, consider the random $K_k$-removal algorithm on $n$ vertices. Let $M$
 be the number of steps it takes the process to terminate and $E(M)$ be the size
 of the resulting $K_k$-free graph. With high probability, $|E(M)|\leqslant n^{2-1/(k(k-1)-2)+o(1)}$.
 \end{theorem}
\noindent Though our bound exists a gap with $|E(M)|\leqslant n^{2k/(k+1)+o(1)}$
conjectured in~\cite{bennett15},
we will show our result corresponds to the inherent barrier of the algorithm.

The remainder of this paper is organized as follows.
In the next section, notations and some lemmas  for
analyzing the random $K_k$-removal algorithm  are presented.
 In Section 3, we discuss the evolution of the algorithm in detail and
estimate the trajectories of these random variables.
 We formally prove the concentrations in Section 4.

\section{Notations and Some Lemmas}

Let $(\Omega,\mathcal{F},\mathbb{P})$ be an arbitrary probability space.
Note that our probability space is the set of
all maximal sequences of edge-disjoint $K_k$s on vertex set $[n]$ with probability
measure given by the uniform random choice at each step. Let
$(\mathcal{F}_i)_{i\geqslant 0}$ be the filtration given by the evolutionary algorithm.
Given a sequence of random variables $X_i$, let $\Delta X=X_{i+1}-X_i$ denote the one-step change 
 for the random variables $X_i$ and the pair $\{X_i,\mathcal{F}_i\}_{i\geqslant 0}$
 is then called a submartingale (resp. supermartingale)
if $X_i$ is $\mathcal{F}_i$-measurable  and $\mathbb{E}[\Delta X|{\mathcal{F}}_i]\geqslant 0$
(resp. $\mathbb{E}[\Delta X|{\mathcal{F}}_i]\leqslant 0$) for all $i\geqslant 0$.
An event is said to occur  with high probability (w.h.p. {\it for short}),
if the probability that it holds tends to 1 when $n\rightarrow\infty$.
Furthermore, for two positive-valued functions $f, g$ on the variable $n$,
we write  $f\ll g$ to denote $\lim_{n\rightarrow \infty}f(n)/g(n)=0$ and $f\sim g$ to denote
$\lim_{n\rightarrow \infty}f(n)/g(n)=1$. Let $a=b\pm c$ be short for $a\in [b-c,b+c]$,
${S\choose b}=\emptyset$ and ${a\choose b}=0$ if $b>|S|$ and $b>a$.
We also use the standard asymptotic notation $o$, $O$, $\Omega$ and $\Theta$.
All logarithms are natural, and the floor and ceiling signs are omitted
whenever they are not crucial. Throughout the following sections we assume that $n\rightarrow\infty$.

For $2\leqslant m\leqslant k$, $u\in[n]$ and $U_m=\{u_1,\cdots,u_m\}\in \binom{[n]}{m}$,
let $N_{u}=N_{u}(i)=\{x\in [n]:xu\in E(i)\}$, $N_{U_m}=N_{U_m}(i)=\cap_{i=1}^{m} N_{u_i}$
and $\mathcal{K}_m(i)$ be the set of complete graph $K_m$ in $G(i)$.
Our goal
is to estimate the number of $K_k$s in $G(i)$, that is $|\mathcal{K}_k(i)|$,
which is particularly denoted by $\mathbf{Q}_k(i)$. Define the random variable $\mathbf{R}_{k,U_m}(i)$
to be
\begin{align}
\mathbf{R}_{k,U_m}(i)=\begin{cases}\bigl|{K}_{k-m}\cap {N_{U_m}\choose k-m}\bigr|,
\ 2\leqslant m\leqslant k-1;\\{\textbf 1}_{U_k},\ {m=k.}\end{cases}
\end{align}
 For $2\leqslant m\leqslant k-1$, $\mathbf{R}_{k,U_m}(i)$ counts the number of $K_{k-m}$s in $G(i)$ such that
every vertex in $K_{k-m}$  is in $N_{U_m}$;  particularly
$\mathbf{R}_{k,U_{k-1}}(i)=|N_{U_{k-1}}|$ is the codegree of the vertex
subset $U_{k-1}$.
$\textbf{1}_{U_k}$ is the indicator random variable  with
$\textbf{1}_{U_k}=1$ if the subgraph induced by $U_k$ in $G(i)$ is
complete, instead $\textbf{1}_{U_k}=0$ otherwise. 
 Bennett et al.~\cite{bennett15} ever added more assumptions on
 codegrees  of larger vertex subsets to obtain stronger results on  random greedy hypergraph
matching algorithm. Sometimes for shorthand we will suppress $i$. These random variables in (2.1)  yield important information
about the underlying process.

Suppose that the vertex set of the $(i+1)$-th taken $K_k$ is denoted by  $U_k$.
Let $U_m\in {U_k\choose m}$  with $2\leqslant m\leqslant k$ and
\begin{align*}
\mathbf{Q}_k^{U_m}(i)=\bigl|\{K_k\in G(i)|K_k\cap U_k=U_m\}\bigr|,
\end{align*}namely, $\mathbf{Q}_k^{U_m}(i)$ denotes the number of $K_k$s in $G(i)$ that exactly contains
the vertices $U_m$ in $U_k$. In particular, $\mathbf{Q}_k^{U_k}(i)=1$. Thus, we have
\begin{align}
\mathbf{Q}_k(i)-\mathbf{Q}_k(i+1)
=\sum_{m=2}^{k}\Bigl(\sum_{U_m\in \binom{U_k}{m}}\mathbf{Q}_k^{U_m}(i)\Bigr).
\end{align}
It is observed that $\mathbf{R}_{k,U_m}$ in (2.1) denotes the number of
extensions to one copy of $K_k$ from $U_m$ when $U_m$ is complete.
By inclusion-exclusion formula, we have
\begin{align}
\mathbf{Q}_k^{U_m}(i)&=\mathbf{R}_{k,U_m}+\sum_{T_1\in \binom{U_k\setminus U_m}{1}}(-1)^1\mathbf{R}_{k,U_m\cup T_1}+\cdots+\notag\\
&\quad\quad\sum_{T_{k-m-1}\in \binom{U_k\setminus U_m}{k-m-1}}
(-1)^{k-m-1}\mathbf{R}_{k,U_m\cup T_{k-m-1}}+(-1)^{k-m}\mathbf{R}_{k,U_k}.
\end{align}

Note that
\begin{align*}
\sum_{U_m\in \binom{U_k}{m}}\Bigl(\sum_{T_i\in \binom{U_k\backslash U_m}{i}}
\mathbf{R}_{k,U_m\cup T_i}\Bigr)
= \binom{m+i}{m}\sum_{U_{m+i}\in \binom{U_k}{m+i} }\mathbf{R}_{k,U_{m+i}}
\end{align*}
  for $0\leqslant i\leqslant k-m$ because  each element $\mathbf{R}_{k,U_{m+i}}$
on the right side is counted $ \binom{m+i}{m}$ times on the left side.
Sum the above corresponding displays (2.3) for all $U_m\in \binom{U_k}{m}$
with $2\leqslant m\leqslant k$ altogether into the equation (2.2), then it follows that
\begin{align*}
&\mathbf{Q}_k(i)-\mathbf{Q}_k(i+1)\\
=&\sum_{U_2\in \binom{U_k}{2}}\mathbf{R}_{k,U_2}+\sum_{U_3\in \binom{U_k}{3}}
\biggl[(-1)^1{3\choose 2}+(-1)^0{3\choose 3}\biggr]\mathbf{R}_{k,U_3}+\cdots\\
&\quad\quad+\sum_{U_{k-1}\in{U_k\choose k-1}}\biggl[(-1)^{k-3}{k-1\choose 2}
+\cdots+(-1)^0{k-1\choose k-1}\biggr]\mathbf{R}_{k,U_{k-1}}\\
&\quad\quad+\biggl[(-1)^{k-2}{k\choose 2}+\cdots+(-1)^0{k\choose k}\biggr]\mathbf{R}_{k,U_k}.
\end{align*}
Since
$\sum_{j=2}^r(-1)^{r-j}{r\choose j}=(-1)^r(r-1)$ for any given integer $r\geqslant2$
and $\mathbf{R}_{k,U_k}=1$  in (2.1),
\begin{align}
\mathbf{Q}_k(i)-\mathbf{Q}_k(i+1)
&=\sum_{U_2\in{U_k\choose 2}}\mathbf{R}_{k,U_2}-2\sum_{U_3\in{U_k\choose 3}}\mathbf{R}_{k,U_3}+\cdots\notag\\
&\quad\quad+(-1)^{k-1}(k-2)\sum_{U_{k-1}\in{U_k\choose k-1}}\mathbf{R}_{k,U_{k-1}}+(-1)^{k}(k-1).
\end{align}
Thus,  the expectation
$\mathbb{E}[\Delta \mathbf{Q}_k| \mathcal{F}_i]$ of $\Delta \mathbf{Q}_k$ is
\begin{align}
&\mathbb{E}[\Delta \mathbf{Q}_k| \mathcal{F}_i]\notag\\
&=-\sum\limits_{U_k\in \mathcal{K}_k(i)} \frac{\sum_{U_2\in{U_k\choose 2}}
\mathbf{R}_{k,U_2}+\cdots+(-1)^{k-1}(k-2)\sum_{U_{k-1}\in{U_k\choose k-1}}
\mathbf{R}_{k,U_{k-1}}+(-1)^{k}(k-1)}{\mathbf{Q}_k(i)}\notag\\
&=(-1)^{k+1}(k-1)-\frac{1}{\mathbf{Q}_k(i)}\sum_{U_2\in \mathcal{K}_2(i)}
\mathbf{R}_{k,U_2}^2+\cdots+\frac{(-1)^k(k-2)}{\mathbf{Q}_k(i)}
\sum_{U_{k-1}\in \mathcal{K}_{k-1}(i)}\mathbf{R}_{k,U_{k-1}}^2,
\end{align}
where the last equality is true because
\begin{align*}
\sum_{U_k\in \mathcal{K}_k(i)}\sum_{U_m\in \binom{U_k}{m}}
\mathbf{R}_{k,U_m}=\sum_{U_m\in \mathcal{K}_m(i)}\mathbf{R}_{k,U_m}^2
\end{align*}
for $2\leqslant m\leqslant k-1$ by double counting.

We also need the following lemmas to
establish dynamic concentrations on variables $\mathbf{Q}_k(i)$ and $\mathbf{R}_{k,U_m}$
for any $U_m\in \binom{[n]}{m}$
with $2\leqslant m\leqslant k-1$, which were also used in
~\cite{bennett15,bohman09,bohman15,bohman101,bohman19,pico142,war142}.

 \begin{lemma}[Bohman et al.~\cite{bohman15}]
 Let $a_1,\cdots,a_\ell\in \mathbb{R}$ and some
 $a\in \mathbb{R}$. Suppose that $|a_i-a|\leqslant \varepsilon$ for all $1\leqslant i\leqslant \ell$, then
 $ \frac{(\sum_{i=1}^\ell a_i)^2}{\ell} \leqslant
 \sum_{i=1}^\ell a_i^2\leqslant \frac{(\sum_{i=1}^\ell a_i)^2}{\ell}+4\ell\varepsilon^2.$
\end{lemma}

\begin{lemma}[Hoeffding and Azuma~\cite{hoeffding63}]
Suppose a sequence of random variables $\{X_i\}_{i\geqslant 0}$
is a supermartingale (resp. submartingale) and $|X_{i}-X_{i-1}|<c_i$,
then for any positive integer $\ell$ and any positive real number $a$,
$\mathbb{P}\bigl[X_\ell-X_0\geqslant a\bigr]\leqslant\exp\bigl[ \frac{-a^2}{2\sum_{i=1}^{\ell}c_i^2}\bigr].
\bigl({\mbox{resp.}}\ \mathbb{P}\bigl[X_\ell-X_0\leqslant -a\bigr]\leqslant \exp\bigl[ \frac{-a^2}{2\sum_{i=1}^{\ell}c_i^2}\bigr].\bigr)$
\end{lemma}

Let $\eta,N>0$ be constants.  A sequence of random variables $\{X_i\}_{i\geqslant 0}$ is  $(\eta,N)$-bounded
if $X_i-\eta\leqslant X_{i+1}\leqslant X_i+N$ for all $i\geqslant 0$. For $(\eta,N)$-bounded
supermartingales and submartingales, Bohman~\cite{bohman09} showed that
\begin{lemma}[Bohman~\cite{bohman09}]
Suppose $\{X_i\}_{i\geqslant 0}$ is an $(\eta,N)$-bounded
supermartingale (resp. submartingale) with initial value $0$
and $\eta\leqslant \frac{N}{10}$. Then
for any positive integer $\ell$ and any positive real number $a$ with $a<\eta \ell$,
$\mathbb{P}\bigl[X_\ell\geqslant a\bigr]\leqslant\exp\bigl[-\frac{a^2}{3\ell\eta N}\bigr]. \bigl({\mbox{resp.}}\ \mathbb{P}\bigl[X_\ell\leqslant-a\bigr]\leqslant\exp\bigl[-\frac{a^2}{3\ell\eta N}\bigr].\bigr)$
\end{lemma}

Finally, in order to explain it is definitely possible to further improve our results. 
the lemma below 
in~\cite{cher52} is also required.

\begin{lemma}[\cite{cher52}]
For $X\sim \Bin(n, p)$ and any $0<\xi\leqslant np$, $\mathbb{P}\bigl[|X - np|>\xi\bigr]<2\exp\left[-\xi^2/\left(3np\right)\right]$.

\end{lemma}


\section{Estimates on the variables in $G(i)$}

 In the following, we use some heuristics to anticipate the likely values of the auxiliary
 random variables throughout the process.
 We assume the random $K_k$-removal algorithm produces a graph whose variables are roughly
 the same as they would be in a random graph $\mathcal{G}(n,p)$ with the same edge density.
The classical Erd\H{o}s-R\'{e}nyi random graph $\mathcal{G}(n,p)$ is on vertex set $[n]=\{1,\cdots,n\}$
 and any two vertices appear as an edge independently with probability $p$.

In order to describe the expected trajectories of $\mathbf{Q}_k(i)$ and $\mathbf{R}_{k,U_m}$
as smooth functions for any $U_m=\{u_1, u_2,\cdots, u_m\}\in \binom{[n]}{m}$ with
  $2\leqslant m\leqslant k-1$, we  appropriately rescale the number of steps $i$
  to be $t=t(i)= \frac{i}{n^2}$
 and introduce a notion of edge density as
\begin{align}
p=p(i,n)=1- \frac{k(k-1)i}{n^2}=1-k(k-1)t.
\end{align}
Note that $p$ can be viewed as either a continuous function
of $t$ or as a function of the discrete variable $i$.
We pass between these interpretations without comment.
With this notation, we have
\begin{align}
|E(i)|={n\choose 2}-{k\choose 2}i={n\choose 2}- \frac{1}{2}(1-p)n^2= \frac{1}{2}(n^2p-n)
\end{align}
such that the number of edges in $G(i)$ with edge density $p$
is approximately equal to the one in
the Erd\H{o}s-R\'{e}nyi graph $\mathcal{G}(n,p)$
up to the negligible linear term when $p$ lies in some range.

For a fixed integer $k\geqslant 4$,
$2\leqslant m\leqslant k-1$ and $U_m\in \binom{[n]}{m}$, under the assumption that $G(i)$  resembles $\mathcal{G}(n,p)$,
we anticipate that the expressions of $\mathbf{Q}_k(i)$ and $\mathbf{R}_{k,U_m}$
are 
\begin{align*}
\mathbf{Q}_k(i)\sim \frac{n^k}{k!}p^{{k\choose 2}}
\quad {\rm{and}}\quad \mathbf{R}_{k,U_m}
\sim \frac{n^{k-m}}{(k-m)!}p^{{k\choose 2}-{m\choose 2}},
\end{align*}
where $ \frac{n^k}{k!}p^{{k\choose 2}}$ counts the expected
number of $K_k$s in $\mathcal{G}(n,p)$; ${n-m\choose k-m}p^{{k\choose 2}-{m\choose 2}}
\sim \frac{n^{k-m}}{(k-m)!}p^{{k\choose 2}-{m\choose 2}}$ counts the expected
number of $K_{k-m}$s
in which every vertex is in $N_{U_m}$.
Our main theorem is as follows:
\begin{theorem}
Given a fixed integer $k\geqslant 4$,
let $U_m\in {[n]\choose m}$ 
with $2\leqslant m\leqslant k-1$, then
there exist absolute constants  $\mu$,  $\gamma_m$ and $\lambda$
  such that, with high probability,
\begin{align}
\mathbf{Q}_k(i)&\leqslant \frac{n^k}{k!}p^{{k\choose 2}}+ \frac{n^{k-1}}{2}p^{{k\choose 2}-4},\\
\mathbf{Q}_k(i)&\geqslant \frac{n^k}{k!}p^{{k\choose 2}}-\sigma^2n^{\alpha}p^{-1}\log^\mu n,\\
\mathbf{R}_{k,U_m}&= \frac{n^{k-m}}{(k-m)!}p^{{k\choose 2}-{m\choose 2}}
\pm \sigma n^{\beta_m}\log^{\gamma_m} n
\end{align}
holding for every $i\leqslant i_0$ with $i_0=
\frac{n^2}{k(k-1)}- \frac{\sqrt[3]{2}}{k(k-1)}n^{2- \frac{1}{k(k-1)-2}}\log^{\lambda}n$, where
\begin{align}
\alpha&=k- \frac{{k\choose 2}+1}{2{k\choose 2}-2},\\
\beta_m&=k- m-  \frac{\binom{k}{2}-\binom{m}{2}}{2{k\choose 2}-2}, 
\end{align}
and the error function $\sigma=\sigma(t)$ is taken with initial
value $\sigma(0)=1$ that slowly grows to be
\begin{align}
\sigma=\sigma(t)=1- \frac{k(k-1)}{4}\log p(t).
\end{align}
\end{theorem}

Theorem 3.1 is proved in Section 4. It implies that for these specific
choices of constants satisfying the equations in (3.6) and (3.7),
and the error function $\sigma$ in (3.8), these random variables are around the
heuristical trajectories to the stopping time
$\tau=i_0$ with high probability. These dynamic concentrations
in turn show that the algorithm produces a graph of size at most $|E(i_0)|$ with high probability.
We make no attempt to optimize
 the constants $\mu$,  $\lambda$ and $\gamma_m$
in all error terms with $2\leqslant m\leqslant k-1$. There are many choices of
them that can be balanced to satisfy certain inequalities, such as
$[{k\choose 2}+1]\lambda>\mu+2,
[{k\choose 2}- \binom{m}{2}]\lambda>\gamma_m+1$ with $2\leqslant m\leqslant k-1$,
and $\gamma_2> \frac{1}{2}$, can support our analysis of Theorem 3.1.
 We do not replace them  with their actual values. This is for the interest of
understanding the role of these constants played in the calculations.

\begin{proof}[Proof of Theorem 1.1]
We recover the number of edges when $p=p_0$ to be
\begin{align*}
|E(i_0)|={n \choose 2}-{k\choose 2}i_0\sim \frac{\sqrt[3]{2}}{2}n^{2- \frac{1}{k(k-1)-2}}\log^{\lambda}n.
\end{align*}
Theorem 1.1 follows directly from Theorem 3.1 by $|E(M)|\leqslant |E(i_0)|$
with room to sparse in the power of the logarithmic factor.
\end{proof}

\begin{remark}The variation equations in (3.3)-(3.5) are verified in a straightforward
manner below. According to (3.1), define
\begin{align}
p_0&=p(i_0,n)=1-\frac{k(k-1) i_0}{n^2}=\sqrt[3]{2}n^{- \frac{1}{k(k-1)-2}}\log^{\lambda}n.
\end{align}
Since $ i\leqslant i_0$ in Theorem 3.1, we have $ p\geqslant p_0$ in (3.9). Note that
$\frac{n^k}{k!}p^{{k\choose 2}}\gg\sigma^2n^{\alpha}p^{-1}\log^\mu n$
when  $\alpha$ is in (3.6), and
 appropriate choices of $\lambda$  and $\mu$. It follows that
$\mathbf{Q}_k(i)= (1+o(1))n^k p^{{k\choose 2}}/k!$  in (3.3) and (3.4).
Similarly all the error terms in (3.5) are negligible compared to their respective corresponding
main terms.
\end{remark}

\begin{remark}
Our bound in Theorem 1.1 exists a gap with $|E(M)|\leqslant n^{2k/(k+1)+o(1)}$
conjectured in~\cite{bennett15}.
In fact, the term $n^{2-1/(k(k-1)-2)}$ corresponds to a natural barrier
 in the random $K_k$-removal algorithm.
To illustrate this, as stated in Theorem 3.1,
 we know $G(i)$ is roughly the same with $\mathcal{G}(n,p)$,
 while we notice that the standard variations of $\mathbf{R}_{k,U_m}$
for any $U_m\in \binom{[n]}{m}$ with $2\leqslant m\leqslant k-1$ would be as large as their main
trajectories when  $p$ is around $n^{-1/(k(k-1)-2)}$ (up to logarithmic factors),
which means that the control over $\mathbf{R}_{k,U_m}$ for any $U_m\in \binom{[n]}{m}$ is lost.
\end{remark}

\begin{remark}
As stated in Theorem 3.1,
 we know $G(i)$ is roughly the same with $\mathcal{G}(n,p)$ for $i\leqslant i_0$. Thus,
 when $p$ is around $p_0$ in (3.9),
by a union bound,  it follows that the probability that there exists one $U_m\in \binom{[n]}{m}$
with some $m$ satisfying $2\leqslant m\leqslant k-1$ such that
$|\mathbf{R}_{k,U_m}- \frac{n^{k-m}}{(k-m)!}p_0^{{k\choose 2}- \binom{m}{2}}
|> \sigma n^{\beta_m}\log^{\gamma_m} n$ is at most
\begin{align}
&\sum_{U_m\in \binom{[n]}{m},2\leqslant m\leqslant k-1}\mathbb{P}\biggl[\biggl|\mathbf{R}_{k,U_m}- \frac{n^{k-m}}{(k-m)!}p_0^{{k\choose 2}- \binom{m}{2}}
\biggr|> \sigma n^{\beta_m}\log^{\gamma_m} n\biggr]\notag\\
&<2 \sum_{m=2}^{k-1}\binom{n}{m}\exp\biggl[-\frac{(k-m)!(\sigma n^{\beta_m}\log^{\gamma_m} n)^2}{3n^{k-m}p_0^{{k\choose 2}- \binom{m}{2}}}\biggr]\notag\\
&=\sum_{m=2}^{k-1}\binom{n}{m}\exp\biggl[-\Theta\biggl(n^{k-m- \frac{{k\choose 2}- \binom{m}{2}}{2{k\choose 2}-2}}\biggr)\biggr]
\end{align}
by applying Lemma 2.4 with $\xi=\sigma n^{\beta_m}\log^{\gamma_m} n$,
where the last equality is true because  $\beta_m$ is in (3.7). Since the summand in (3.10) is increasing
in $m$ for fixed $k\geqslant 4$, it suffices to take the number of terms times the last term when $m=k-1$. Thus, we have
\begin{align*}
&\sum_{m=2}^{k-1}\binom{n}{m}\exp\biggl[-\Theta\biggl(n^{k-m- \frac{{k\choose 2}- \binom{m}{2}}{2{k\choose 2}-2}}\biggr)\biggr]\\
&=O\biggl(n^{k-1}\exp\biggl[-\Theta\biggl(n^{1- \frac{{k\choose 2}- \binom{k-1}{2}}{2{k\choose 2}-2}}\biggr)\biggr]\biggr)=o(1).
\end{align*}
In fact, we could show the similar phenomenon even when we  take $\xi=\Theta(n^\theta)$
with $\frac{1}{2}\beta_m <\theta<\beta_m$, 
instead our main results in Theorem 3.1 cannot support us.
Like~\cite{bohman15}, in order to prove better
 bounds on $|E(M)|$, it is possible  to design new 
random variables such that their variations decrease as the process evolves.
\end{remark}

\section{Proof of Theorem 3.1}

Recall the outline of the critical interval method~\cite{bohman101,bohman15,bennett15}
to  control some graph parameters when the process evolves.
Let the stopping time $\tau$ be the minimum of
$i_0$ and the smallest index $i$ such that any one of the random variables
violates its corresponding  trajectory.
Let the event $\mathcal{E}_X$ be of the form $X(i)= x(i)\pm e(i)$ for all $i\leqslant i_0$,
where $X(i)$ is some random variable, $x(i)$ is the expected trajectory and $e(i)$ is the error term.
We show that the event $\{\tau=i_0\}$ holds by means of
$\{\tau=i_0\}=\cap_{X\in\mathcal{I}}\mathcal{E}_X$,
where $|\mathcal{I}|$ is polynomial in $n$.

For each such random variable $X(i)$,
we define a critical interval $I_{X}$ for its bound (upper and lower)
that has one endpoint at the bound we are trying to maintain
and the other slightly closer to the expected trajectory of the random variable.
Consider a fixed step $j\leqslant i_0$ such that $X(j)\in I_{X}$.
Define the stopping time $\tau_{X,j}$ to be
$\tau_{X,j}=\min\{i_0,\max\{j,\tau\},{\mbox{the\ smallest\ }}i\geqslant
j\ {\mbox{ such\  that}\ }X(i)\notin I_{X}\}$,
which  made us possible to establish the martingale condition
and apply the martingale inequality in Lemma 2.2 or Lemma 2.3.
Establish bounds on the events that the designed variable
 crosses its critical interval in the process,
such that a simple application of the union bound over all starting point $j$ 
shows that the probability of the occurrence of any event in the collection is low
to complete the proof.

As a supplement,  we list some necessary
inequalities that we need in the following proof of  Theorem 3.1.
By Lemma 2.1, we have
\begin{align*}
\sum_{U_m\in \mathcal{K}_m(i)}\mathbf{R}_{k,U_m}^2\geqslant \frac{
(\sum_{U_m\in \mathcal{K}_m(i)}\mathbf{R}_{k,U_m})^2}{
|\mathcal{K}_m(i)|}
\end{align*} for any $U_m\in \mathcal{K}_m(i)$ with $2\leqslant m\leqslant k-1$.
Firstly, note that $\sum_{U_m\in \mathcal{K}_m(i)}\mathbf{R}_{k,U_m}={k\choose m}\mathbf{Q}_k(i)$
because each element on the right side is counted ${k \choose m}$ times on the left side.
Next, note that $|\mathcal{K}_{2}(i)|=|E(i)|\sim \frac{n^2}{2}p$ in (3.2) when $p\geqslant p_0$ in (3.9), and
we recursively apply the equation $|\mathcal{K}_{m}(i)|\leqslant \frac{n}{m}|\mathcal{K}_{m-1}(i)|$
to achieve $|\mathcal{K}_{m}(i)|\leqslant\frac{n^m}{m!}p$ with $2\leqslant m\leqslant k-1$.
Thus, we have
\begin{align}
\sum_{U_m\in \mathcal{K}_m(i)}\mathbf{R}_{k,U_m}^2\geqslant
\frac{m!{k\choose m}^2\mathbf{Q}_k^2(i)}{n^mp}.
\end{align}

Conditioned on the estimates in (3.5) hold on $\mathbf{R}_{k,U_m}$ for
any $U_m\in {[n]\choose m}$ with $2\leqslant m\leqslant k-1$, we also have the upper
bounds of $\sum_{U_m\in \mathcal{K}_m(i)}\mathbf{R}_{k,U_m}^2$.
For $m=2$, we have $\beta_2=k- \frac{5}{2}$ in (3.7) and $|\mathcal{K}_{2}(i)|
\sim \frac{n^2}{2}p$, then by Lemma 2.1,
\begin{align}
\sum_{U_2\in \mathcal{K}_2(i)}\mathbf{R}_{k,U_2}^2&\leqslant
\frac{\bigl(\sum_{U_2\in \mathcal{K}_2(i)}\mathbf{R}_{k,U_2}\bigr)^2}{|\mathcal{K}_2(i)|}+
4|\mathcal{K}_2(i)|\bigl(\sigma n^{\beta_2}\log^{\gamma_2} n\bigr)^2\notag\\
&\sim \frac{2!{k\choose 2}^2\mathbf{Q}_k^2(i)}{n^2p}+2\sigma^2 n^{2k-3}p\log^{2\gamma_2} n.
\end{align}
For $3\leqslant m\leqslant k-1$, 
by the estimates in (3.5) and $|\mathcal{K}_{m}(i)|\leqslant\frac{n^m}{m!}p$, the trivial upper bound is
\begin{align}
\sum_{U_m\in \mathcal{K}_m(i)}\mathbf{R}_{k,U_m}^2\leqslant
\frac{n^mp}{m!}\Bigl(\frac{n^{k-m}}{(k-m)!}p^{{k\choose 2}-{m\choose 2}}+\sigma n^{\beta_m}\log^{\gamma_m} n\Bigr)^2.
\end{align}

\subsection{Tracking  $\mathbf{Q}_k(i)$}

For the upper bound of $\mathbf{Q}_k(i)$, we introduce a critical interval as
\begin{align}
I_{\mathbf{Q}_k}^u=\Bigl( \frac{n^k}{k!}p^{k\choose 2}+Bn^{k-1}p^{{k\choose 2}-4},
\frac{n^k}{k!}p^{k\choose 2}+ \frac{n^{k-1}}{2}p^{{k\choose 2}-4}\Bigr),
\end{align}
where
\begin{align}
B&= \frac{1}{2}- \frac{1}{2{k\choose 2}}+ \frac{1}{3{k\choose 2}(k-4)!}<\frac{1}{2}.%
\end{align}
Consider a fixed step $j\leqslant i_0$. Suppose $\mathbf{Q}_k(j)\in I_{\mathbf{Q}_k}^u$. Define
\begin{align}
\tau_{\mathbf{Q}_k,j}^u=\min\bigl\{i_0,\max\{j,\tau\},{\mbox{the\ smallest\ }}i\geqslant
j{\mbox{ such\  that}\ }\mathbf{Q}_k(i)\notin I_{\mathbf{Q}_k}^u\bigr\}.
\end{align}
Let $j\leqslant i\leqslant\tau_{\mathbf{Q}_k,j}^u$, thus  all calculations in this subsection
are conditioned on the estimates in (3.5) hold on $\mathbf{R}_{k,U_m}$
for any $U_m\in {[n]\choose m}$ with $2\leqslant m\leqslant k-1$.

By the equation shown in (2.5),  it follows that
\begin{equation*}
\begin{aligned}
\mathbb{E}[\Delta \mathbf{Q}_k| \mathcal{F}_i]
&=(-1)^{k+1}(k-1)- \frac{1}{\mathbf{Q}_k(i)}\sum_{U_2\in  \mathcal{K}_2(i)}
\mathbf{R}_{k,U_2}^2+\cdots+
\frac{(-1)^k(k-2)}{\mathbf{Q}_k(i)}\sum_{U_{k-1}\in \mathcal{K}_{k-1}(i)}\mathbf{R}_{k,U_{k-1}}^2\\
&<(-1)^{k+1}(k-1)- \frac{2{k\choose 2}^2\mathbf{Q}_k(i)}{n^2p}+ \frac{2}{\mathbf{Q}_k(i)}\frac{n^3 p}{3!}
\Bigl( \frac{n^{k-3}}{(k-3)!}p^{{k\choose 2}-3}+\sigma n^{\beta_3}\log^{\gamma_3} n\Bigr)^2\\
&\qquad+O\bigl(n^{k-4}p^{{k\choose 2}-1}\bigr),
\end{aligned}
\end{equation*}
where  $\sum_{U_2\in  \mathcal{K}_2(i)}
\mathbf{R}_{k,U_2}^2$ and $\sum_{U_3\in  \mathcal{K}_3(i)}
\mathbf{R}_{k,U_3}^2$ are replaced by the equations in (4.1) and (4.3), the last term $O(n^{k-4}p^{{k\choose 2}-1})$ comes from
$\sum_{U_4\in \mathcal{K}_4(i)}\mathbf{R}_{k,U_4}^2$ in (4.1)
that dominates all the remaining terms.

Since $\mathbf{Q}_k(i)\in I_{\mathbf{Q}_k}^u$  is in (4.4), we further have
\begin{align}
\mathbb{E}[\Delta \mathbf{Q}_k| \mathcal{F}_i]
&<(-1)^{k+1}(k-1)-\frac{2{k\choose 2}^2 n^{k-2}}{k!}p^{{k\choose 2}-1}-2{k\choose 2}^2Bn^{k-3}p^{{k\choose 2}-5}\notag\\
&\quad\quad+ \frac{k!n^{k-3}}{3(k-3)!^2}p^{{k\choose 2}-5}+O\bigl(\sigma n^{\beta_3}p^{-2}\log^{\gamma_3}n\bigr),
\end{align}
where $O(n^{k-4}p^{{k\choose 2}-1})$ is absorbed into $O(\sigma n^{\beta_3}p^{-2}\log^{\gamma_3}n)$
when $\beta_3$ is in (3.7).

For all $i$ with $j\leqslant i\leqslant\tau_{\mathbf{Q}_k,j}^u$, define the sequence of random variables to be
\begin{align}
\mathbf{U}(i)= \mathbf{Q}_k(i)- \frac{n^k}{k!}p^{k\choose 2}- \frac{n^{k-1}}{2}p^{{k\choose 2}-4}.
\end{align}

\noindent{\bf Claim 4.1:}\  The
sequence $\mathbf{U}(j),\mathbf{U}(j+1),\cdots,\mathbf{U}(\tau_{\mathbf{Q}_k,j}^u)$
is a supermartingale and the maximum one step $\Delta \mathbf{U}$ is
$O(\sigma n^{k- 5/2}\log^{\gamma_2} n)$.

\begin{proof}[Proof of Claim 4.1]\ To see this, for $j\leqslant i\leqslant\tau_{\mathbf{Q}_k,j}^u$,
as the equation in (4.8), we have
\begin{equation*}
\begin{aligned}
\mathbb{E}[\Delta \mathbf{U}|\mathcal{F}_i]
&=\mathbb{E}[\Delta \mathbf{Q}_k|\mathcal{F}_i]- \frac{n^k}{k!}\Bigl[p^{k\choose 2}(i+1)-p^{k\choose 2}(i)\Bigr]-
 \frac{n^{k-1}}{2}\Bigl[p^{{k\choose 2}-4}(i+1)-p^{{k\choose 2}-4}(i)\Bigr].
\end{aligned}
\end{equation*}
Note that $p=p(i)=1-k(k-1)t$, $p(i+1)=1-k(k-1)(t+\frac{1}{n^2})$ in (3.1), then by
Taylor's expansion, we have
\begin{align}
\mathbb{E}[\Delta \mathbf{U}|\mathcal{F}_i]
&=\mathbb{E}[\Delta \mathbf{Q}_k|\mathcal{F}_i]- \frac{n^k}{k!}\biggl[-{k\choose 2}
\frac{k(k-1)}{n^2}p^{{k\choose 2}-1}+O\Bigl( \frac{1}{n^4}p^{{k\choose 2}-2}\Bigr)\biggr]\notag\\
&\quad\quad- \frac{n^{k-1}}{2}\biggl[-\biggl({k\choose 2}-4\biggr) \frac{k(k-1)}{n^2}p^{{k\choose 2}-5}
+O\Bigl(\frac{1}{n^4}p^{{k\choose 2}-6}\Bigr)\biggr]\notag\\
&=\mathbb{E}[\Delta \mathbf{Q}_k|\mathcal{F}_i]+ \frac{2{k\choose 2}^2 n^{k-2}}{k!}p^{{k\choose 2}-1}+
\biggl[{k\choose 2}-4\biggr]{k\choose 2}n^{k-3}p^{{k\choose 2}-5}\notag\\
&\quad\quad+O\bigl(n^{k-4}p^{{k\choose 2}-2}\bigr),
\end{align}
where  $O(n^{k-5}p^{{k\choose 2}-6})$ is absorbed into
$O(n^{k-4}p^{{k\choose 2}-2})$ when $p\geqslant p_0$ in (3.9). With the help of the equation in (4.7), we further have
\begin{equation*}
\begin{aligned}
\mathbb{E}[\Delta \mathbf{U}|\mathcal{F}_i]
&<(-1)^{k+1}(k-1)-\biggl[2{k\choose 2}^2B-{k\choose 2}^2+4{k\choose 2}-
\frac{k!}{3(k-3)!^2}\biggr]n^{k-3}p^{{k\choose 2}-5}\\
&\quad\quad+O\bigl( \sigma n^{\beta_3}p^{-2}\log^{\gamma_3}n\bigr)\\
&<(-1)^{k+1}(k-1)-2{k\choose 2}n^{k-3}p^{{k\choose 2}-5}+
O\bigl( \sigma n^{\beta_3}p^{-2}\log^{\gamma_3}n\bigr),
\end{aligned}
\end{equation*}
where
\begin{align*}2{k\choose 2}^2B-{k\choose 2}^2+4{k\choose 2}- \frac{k!}{3(k-3)!^2}=3{k\choose 2}
-{k\choose 2} \frac{2}{3(k-3)!}>2{k\choose 2}
\end{align*}  by
$B$ shown in (4.5), and $O(n^{k-4}p^{{k\choose 2}-2})$ is absorbed into
$O(\sigma n^{\beta_3}p^{-2}\log^{\gamma_3}n)$ by
$\beta_3$ shown in (3.7). Note that ${k\choose 2}n^{k-3}p^{{k\choose 2}-5}>
O( \sigma n^{\beta_3}p^{-2}\log^{\gamma_3}n)+(-1)^{k+1}(k-1)$ when $p\geqslant p_0$  in (3.9),
and appropriate choices of $\lambda$ and $\gamma_3$,
then we have $\mathbb{E}[\Delta \mathbf{U}|\mathcal{F}_i]<0$ and the
sequence $\mathbf{U}(j),\mathbf{U}(j+1),\cdots,\mathbf{U}(\tau_{\mathbf{Q}_k,j}^u)$ is
a supermartingale.

Next, we show the maximum one step $\Delta \mathbf{U}$ is  $O(\sigma n^{k-5/2}\log^{\gamma_2} n)$.
As the equations shown in (4.8) and  (4.9), we have
\begin{equation*}
\begin{aligned}
\Delta \mathbf{U}&=\Delta \mathbf{Q}_k+ \frac{2{k\choose 2}^2}{k!}n^{k-2}p^{{k\choose 2}-1}+
\biggl[{k\choose 2}-4\biggr]{k\choose 2}n^{k-3}p^{{k\choose 2}-5}
+O\bigl(n^{k-4}p^{{k\choose 2}-2}\bigr).\\
\end{aligned}
\end{equation*}%
Apply the equation of $\Delta \mathbf{Q}_k$ shown in (2.4) to the above display, by
the equation of $\textbf{R}_{k,U_m}$ shown in (3.5) for any $U_m\in {[n]\choose m}$,
and $\beta_m$ shown
in (3.7) with $2\leqslant m\leqslant k-1$,  then  we finally have
\begin{equation*}
\begin{aligned}
\Delta \mathbf{U}&\leqslant-{k\choose 2}\Bigl( \frac{n^{k-2}}{(k-2)!}p^{{k\choose 2}-1}-
\sigma n^{\beta_2}\log^{\gamma_2} n\Bigr)+{k\choose 3}\Bigl( \frac{n^{k-3}}{(k-3)!}
p^{{k\choose 2}-{3\choose 2}}+\sigma n^{\beta_3}\log^{\gamma_3} n\Bigr)+\cdots\\
&\quad\quad+ \frac{2{k\choose 2}^2}{k!}n^{k-2}p^{{k\choose 2}-1}+
\biggl[{k\choose 2}-4\biggr]{k\choose 2}n^{k-3}p^{{k\choose 2}-5}+O\bigl(n^{k-4}p^{{k\choose 2}-2}\bigr)\\
&=O\bigl(\sigma n^{k-\frac{5}{2}}\log^{\gamma_2} n\bigr).
\end{aligned}
\end{equation*}
The claim follows.
 \end{proof}

Now, apply Lemma 2.2 to the sequence $\mathbf{U}(j),\mathbf{U}(j+1),\cdots,
\mathbf{U}(\tau_{\mathbf{Q}_k,j}^u)$. The number of steps in this sequence
 is  $O(n^2p)$ because $|E(i)|\sim \frac{n^2}{2}p$  in (3.2) when $p\geqslant p_0$ in (3.9).
Since $\mathbf{Q}_k(j)\in I_{\mathbf{Q}_k}^u$ in (4.4), we have the initial value
$\mathbf{U}(j)\geqslant -\bigl( \frac{1}{2{k\choose 2}}-
\frac{1}{3{k\choose 2}(k-4)!}\bigr)n^{k-1}p^{{k\choose 2}-4}$.
Then, for all $i$ with $j\leqslant i\leqslant \tau_{\mathbf{Q}_k,j}^u$, the probability of  a large deviation for $\mathbf{Q}_k(i)$
beginning at the step $j$ is at most
\begin{equation*}
\begin{aligned}
&\mathbb{P}\Bigl[\mathbf{Q}_k(i)\geqslant \frac{n^k}{k!}p^{k\choose 2}+
\frac{n^{k-1}}{2}p^{{k\choose 2}-4}\Bigr]\\
&=\mathbb{P}\Bigl[\mathbf{U}(i)\geqslant 0\Bigr]\\
&\leqslant \exp\biggl[-\Omega\biggl( \frac{(n^{k-1}p^{{k\choose 2}-4})^2}{(n^2p)
\bigl(\sigma n^{k-5/2}\log^{\gamma_2} n\bigr)^2}\biggr)\biggr]\\
&=\exp\biggl[-\Omega\biggl( \frac{n p^{2{k\choose 2}-9}}{\sigma^2\log^{2\gamma_2}n}\biggr)\biggr].
\end{aligned}
\end{equation*}
By the union bound, note that there are at most $n^2$ possible values of
$j$ in (3.1) and $p\geqslant p_0$ in (3.9), then we have
\begin{equation*}
\begin{aligned}
n^2\exp\biggl[-\Omega\biggl( \frac{n p^{2{k\choose 2}-9}}{\sigma^2\log^{2\gamma_2}n}\biggr)\biggr]=o(1).
\end{aligned}
\end{equation*}
W.h.p., $\mathbf{Q}_k(i)$ never crosses its critical interval $I_{\mathbf{Q}_k}^u$ in
(4.1), and so the upper bound of $\mathbf{Q}_k(i)$ in (3.3) is true.

\begin{remark}
Proving the lower bound of $\mathbf{Q}_k(i)$  is similar.
 We show the proof in the appendix for reference.
\end{remark}

\subsection{Tracking $\mathbf{R}_{k,U_m}$ for any $U_m\in{[n]\choose m}$ with $2\leqslant m\leqslant k-1$}

We prove the dynamic concentration of  $\mathbf{R}_{k,U_m}$
for any $U_m\in{[n]\choose m}$ with $2\leqslant m\leqslant k-1$ in this subsection.
Fix one subset $U_{m^*}\in{[n]\choose m^*}$ for some $m^*$ with
$2\leqslant m^*\leqslant k-1$.  We start with the upper bound of
$\mathbf{R}_{k,U_{m^*}}$. Our critical interval for the upper bound of $\mathbf{R}_{k,U_{m^*}}$ is
\begin{align}
I_{\mathbf{R}_{k,U_{m^*}}}^u=\Bigl( \frac{n^{k-m^*}}{(k-m^*)!}&p^{{k\choose 2}
-{m^*\choose 2}}+(\sigma-1) n^{\beta_{m^*}}\log^{\gamma_{m^*}} n,\notag\\&\qquad
\frac{n^{k-m^*}}{(k-m^*)!}p^{{k\choose 2}-{m^*\choose 2}}+\sigma n^{\beta_{m^*}}\log^{\gamma_{m^*}} n\Bigr),
\end{align}
where $\beta_{m^*}=k-m^*- \frac{\binom{k}{2}-\binom{m^*}{2}}{2\binom{k}{2}-2}$ in (3.7).
Consider a fixed step $j\leqslant i_0$. Suppose $\mathbf{R}_{k,U_{m^*}}(j)\in I_{\mathbf{R}_{k,U_{m^*}}}^u$. Define
\begin{align}
\tau_{\mathbf{R}_{k,U_{m^*}},j}^u=\min\bigl\{i_0,\max\{j,\tau\},
{\mbox{the\ smallest\ }}i\geqslant j\ {\mbox{ such\  that}\ }
\mathbf{R}_{k,U_{m^*}}\notin I_{\mathbf{R}_{k,U_{m^*}}}^u\bigr\}.
\end{align}
Let $j\leqslant i\leqslant\tau_{\mathbf{R}_{k,U_{m^*}},j}^u$, thus all calculations
are conditioned on the events that the estimates in (3.3) and (3.4) hold on $\mathbf{Q}_k(i)$, and
the estimates in (3.5) hold on $\mathbf{R}_{k,U_m}$
for all  $U_m\in {[n]\choose m}$ with $2\leqslant m\leqslant k-1$ and $U_m\neq U_{m^*}$.

Take one $U_{m^*}^c\in K_{k-m^*}\cap N_{U_{m^*}}$ in $G(i)$
and let $\mathbf{Q}_{k,U_{m^*},U_{m^*}^c}$ be the number of
$K_k$s in $G(i)$ such that the removal of the edges in any one of these $K_k$s
results in $U_{m^*}^c\notin K_{k-m^*}\cap N_{U_{m^*}}$ in $G(i+1)$. Then,
we have
\begin{align}
\mathbb{E}[\Delta \mathbf{R}_{k,{U_{m^*}}}| \mathcal{F}_i]
&=-\sum_{{U_{m^*}^c}\in K_{k-m^*}\cap N_{U_{m^*}}}
\frac{\mathbf{Q}_{k,U_{m^*},U_{m^*}^c}}{\mathbf{Q}_k(i)}.
\end{align}

In order to count $\mathbf{Q}_{k,U_{m^*},U_{m^*}^c}$,
let $H\subseteq U_{m^*}\cup U_{m^*}^c$ and $\mathbf{Q}_{k,U_{m^*},U_{m^*}^c}^{H}$
be the number of $K_k$s in $\mathbf{Q}_{k,U_{m^*},U_{m^*}^c}$
such that these $K_k$s satisfy $K_k\cap (U_{m^*}\cup U_{m^*}^c)=H$.
Define $|H|=h$. To ensure that the removal of the edges in any one of
 these $K_k$s results in $U_{m^*}^c\notin K_{k-m^*}\cap N_{U_{m^*}}$ in $G(i+1)$,
 it is observed that $H\cap U_{m^*}^c\neq \emptyset$ and
$2\leqslant h\leqslant k$.

Choose $H\in {\cup_{\rho=0}^{h-1}{U_{m^*}\choose \rho}\oplus{U_{m^*}^c\choose h-\rho}}$,
where ${U_{m^*}\choose \rho}\oplus{U_{m^*}^c\choose h-\rho}$ denotes the collection of union sets consisting of
$\rho$ vertices in $U_{m^*}$ and $h-\rho$ vertices in $U_{m^*}^c$. Hence,
$\mathbf{Q}_{k,U_{m^*},U_{m^*}^c}$ is decomposed into
\begin{align}
\mathbf{Q}_{k,U_{m^*},U_{m^*}^c}=\sum_{h=2}^{k}\sum_{\rho=0}^{h-1}
\sum_{H\in {{U_{m*}\choose \rho}\oplus{U_{m^*}^c\choose h-\rho}}}\mathbf{Q}_{k,U_{m^*},U_{m^*}^c}^{H}.
\end{align}
Following the  inclusion-exclusion counting technique shown in (2.4), we have
\begin{align*}
\mathbf{Q}_{k,U_{m^*},U_{m^*}^c}^{H}
&=\mathbf{1}_H\cdot\mathbf{R}_{k,H}-\sum_{T_1\in{(U_{m^*}\cup U_{m^*}^c)
\setminus H\choose 1}}\mathbf{1}_{H\cup T_1}\cdot\mathbf{R}_{k,H\cup T_1}+\cdots\\&\quad\quad+
\sum_{T_{k-h}\in{(U_{m^*}\cup U_{m^*}^c)\setminus H\choose k-h}}(-1)^{k-h}\mathbf{1}_{H\cup T_{k-h}}\cdot\mathbf{R}_{k,H\cup T_{k-h}}\\
&=\sum_{z=0}^{k-h}\sum_{T_z\in{(U_{m^*}\cup U_{m^*}^c)\setminus H\choose z}}(-1)^z\mathbf{1}_{H\cup T_z}\cdot\mathbf{R}_{k,H\cup T_z},
\end{align*}
where $\mathbf{1}_{H\cup T_z}$ with $0\leqslant z\leqslant k-h$ is the indicator random variable
 depending on whether the subgraph induced by $H\cup T_z$ in $G(i)$ is complete or not.
 Combining with the equation in (4.13), we further have
\begin{align*}
\mathbf{Q}_{k,U_{m^*},U_{m^*}^c}
&=\sum_{h=2}^{k}\sum_{\rho=0}^{h-1}\sum_{z=0}^{k-h}\sum_{H\in {{U_{m*}\choose \rho}
\oplus{U_{m^*}^c\choose h-\rho}}}\sum_{T_z\in{(U_{m^*}\cup U_{m^*}^c)
\setminus H\choose z}}(-1)^z\mathbf{1}_{H\cup T_z}\cdot\mathbf{R}_{k,H\cup T_z}.
\end{align*}
In the above display, for  fixed integers $h$ and $z$, we
recount the union $H\cup T_z$
as a subset $H_{h+z}\in {U_{m^*}\choose \zeta}
\oplus{U_{m^*}^c\choose h+z-\zeta}$ with $0\leqslant \zeta\leqslant h+z$,
then each $H_{h+z}$ is counted $[{h+z\choose h}-{\zeta\choose h}]$ times in
$H\cup T_z$ because $H\cap U_{m^*}^c\neq \emptyset$, which means that
\begin{align*}
&\sum_{\rho=0}^{h-1}\sum_{H\in {{U_{m*}\choose \rho}
\oplus{U_{m^*}^c\choose h-\rho}}}\sum_{T_z\in{(U_{m^*}\cup U_{m^*}^c)
\setminus H\choose z}}(-1)^z\mathbf{1}_{H\cup T_z}\cdot\mathbf{R}_{k,H\cup T_z}\\
&=\sum_{\zeta=0}^{h+z}\sum_{H_{h+z}\in {U_{m^*}\choose \zeta}
\oplus{U_{m^*}^c\choose h+z-\zeta}}\biggl[{h+z\choose h}-{\zeta\choose h}\biggr]
(-1)^z\mathbf{1}_{H_{h+z}}\cdot\mathbf{R}_{k,H_{h+z}}.
\end{align*}
 It follows that
\begin{align}
\mathbf{Q}_{k,U_{m^*},U_{m^*}^c}
&=\sum_{h=2}^{k}\sum_{z=0}^{k-h}\sum_{\zeta=0}^{h+z}\sum_{H_{h+z}\in {U_{m^*}\choose \zeta}
\oplus{U_{m^*}^c\choose h+z-\zeta}}\biggl[{h+z\choose h}-{\zeta\choose h}\biggr]
(-1)^z\mathbf{1}_{H_{h+z}}\cdot\mathbf{R}_{k,H_{h+z}}.
\end{align}

In fact,  $\mathbf{Q}_{k,U_{m^*},U_{m^*}^c}$ is the sum of all elements
in the upper triangular matrix below 
$$\small{
\left(\begin{array}{ccc}
\sum\limits_{\zeta=0}^{2}\sum\limits_{H_{2}\in {U_{m^*}\choose \zeta}\oplus{U_{m^*}^c\choose 2-\zeta}}(-1)^0[{2\choose 2}-{\zeta\choose 2}]\mathbf{1}_{H_2}\cdot\mathbf{R}_{k,H_{2}}&\cdots&\sum\limits_{\zeta=0}^{k}\sum\limits_{H_{k}\in {U_{m^*}\choose \zeta}\oplus{U_{m^*}^c\choose k-\zeta}}(-1)^{k-2}[{k\choose 2}-{\zeta\choose 2}]\mathbf{1}_{H_k}\cdot\mathbf{R}_{k,H_{k}}\\
\cdots&\cdots
&\cdots\\
\sum\limits_{\zeta=0}^{k}\sum\limits_{H_{k}\in {U_{m^*}\choose \zeta}\oplus{U_{m^*}^c\choose k-\zeta}}(-1)^{0}[{k\choose k}-{\zeta\choose k}]\mathbf{1}_{H_k}\cdot\mathbf{R}_{k,H_{k}}&\cdots&0
\end{array}\right)}
$$
with the line corresponds
to the index $h$ and the column
corresponds to the index $z$ in (4.14), respectively.
Recalculate $\mathbf{Q}_{k,U_{m^*},U_{m^*}^c}$ according to every back diagonal
lines to be 
\begin{align}
\mathbf{Q}_{k,U_{m^*},U_{m^*}^c}&=\sum_{h=2}^{k}\sum_{\zeta=0}^{h}\sum_{H_{h}\in {U_{m^*}\choose \zeta}\oplus{U_{m^*}^c\choose h-\zeta}}\sum_{s=2}^{h}(-1)^{h-i}\biggl[{h\choose s}-{\zeta\choose s}\biggr]\mathbf{1}_{H_h}\cdot\mathbf{R}_{k,H_{h}}.
\end{align}
Note that there is no $\mathbf{R}_{k,U_{m^*}}$ on the right side of (4.15)
because $\mathbf{R}_{k,U_{m^*}}$ corresponds to the case when
$\zeta=h$. Thus, the estimates on $\mathbf{Q}_k(i)$ in (3.3) and (3.4),
the estimates on $\mathbf{R}_{k,U_m}$ in (3.5) for all $U_m\in {[n]\choose m}$,
$U_m\neq U_{m*}$ with $2\leqslant m\leqslant k-1$, already support the calculation
of $\mathbf{Q}_{k,U_{m^*},U_{m^*}^c}$ in (4.15).

Furthermore, according to the expressions of $\mathbf{R}_{k,H_h}$ for $2\leqslant h\leqslant k-1$
in (3.5), the term $\mathbf{R}_{k,H_2}$ dominates the sum on the right side of (4.15).
Thus, we have $\zeta=0,1$ and $s=2$. It follows that, 
\begin{align}
&\mathbf{Q}_{k,U_{m^*},U_{m^*}^c}\notag\\
&=\biggl[{k-m^*\choose 2}+m^*(k-m^*)\biggl]\Bigl( \frac{n^{k-2}}{(k-2)!}p^{{k\choose 2}-1}
- \sigma n^{\beta_2}\log^{\gamma_2} n\Bigr)+O\bigl(n^{k-3}p^{{k\choose 2}-3}\bigr),
\end{align}
where ${k-m^*\choose 2}$ counts the number of $\mathbf{R}_{k,H_{2}}$  when $\zeta=0$ and $s=2$,
$m^*(k-m^*)$ counts the number of $\mathbf{R}_{k,H_{2}}$  when $\zeta=1$ and $s=2$.
Note that ${k-m^*\choose 2}+m^*(k-m^*)={k\choose 2}-{m^*\choose 2}$
and $\beta_2=k- \frac{5}{2}$ in (3.7),
combining the equations in (4.12) and (4.16), and applying the estimates
of $\mathbf{Q}_k(i)$  in (3.3),
we have
\begin{equation*}
\begin{aligned}
\mathbb{E}[\Delta \mathbf{R}_{k,U_{m^*}}| \mathcal{F}_i]
<-\sum_{U_{m^*}^c\in K_{k-m^*}\cap N_{U_{m^*}}} \frac{\bigl[{k\choose 2}-
{m^*\choose 2}\bigl]\bigl( \frac{n^{k-2}}{(k-2)!}p^{{k\choose 2}-1}-
\sigma n^{k- \frac{5}{2}}\log^{\gamma_2} n\bigr)+O\bigl(n^{k-3}p^{{k\choose 2}-3}\bigr)}{\frac{n^k}{k!}p^{{k\choose 2}}}.
\end{aligned}
\end{equation*}
The ways to choose $U_{m^*}^c\in K_{k-m^*}\cap N_{U_{m^*}}$ is $\mathbf{R}_{k,U_{m^*}}$
and $\mathbf{R}_{k,U_{m^*}}\in I_{\mathbf{R}_{k,U_{m^*}}}^u$ in (4.10), then it further follows that
\begin{equation*}
\begin{aligned}
&\mathbb{E}[\Delta \mathbf{R}_{k,U_{m^*}}| \mathcal{F}_i]\\
&<- \frac{\bigl[{k\choose 2}-{m^*\choose 2}\bigr]
\bigl( \frac{n^{k-m^*}}{(k-m^*)!}p^{{k\choose 2}-{m^*\choose 2}}
+(\sigma-1) n^{\beta_{m^*}}\log^{\gamma_{m^*}} n\bigr)
\bigl( \frac{n^{k-2}}{(k-2)!}p^{{k\choose 2}-1}- \sigma n^{k-
\frac{5}{2}}\log^{\gamma_2} n\bigr)}{ \frac{n^k}{k!}p^{{k\choose 2}}}\\
&\quad +O\bigl(n^{k-m^*-3}p^{{k\choose 2}-{m^*\choose 2}-3}\bigr).
\end{aligned}
\end{equation*}
Rearrange the above equation to be
\begin{align}
\mathbb{E}[\Delta \mathbf{R}_{k,U_{m^*}}| \mathcal{F}_i]
&<- \frac{\bigl[{k\choose 2}-{m^*\choose 2}\bigr] k(k-1)n^{k-m^*-2}}{(k-m^*)!}
p^{{k\choose 2}-{m^*\choose 2}-1}\notag\\
&\quad\quad+
\frac{\bigl[{k\choose 2}-{m^*\choose 2}\bigr] k!\sigma n^{k-{m^*}-
\frac{5}{2}}\log^{\gamma_2} n}{(k-m^*)!p^{{m^*\choose 2}}}\notag\\
&\quad\quad- \frac{\bigl[{k\choose 2}-{m^*\choose 2}\bigr]k(k-1)(\sigma-1)}
{p}n^{\beta_{m^*}-2}\log^{\gamma_{m^*}} n
\notag\\
&\quad\quad+
\frac{\bigl[{k\choose 2}-{m^*\choose 2}\bigr]k!\sigma(\sigma-1)}{p^{{k\choose 2}}}n^{\beta_{m^*}-\frac{5}{2}}
\log^{\gamma_2+\gamma_{m^*}} n\notag\\
&\quad\quad+O\bigl(n^{k-m^*-3}p^{{k\choose 2}-{m^*\choose 2}-3}\bigr).
\end{align}

For all $i$ with $j\leqslant i\leqslant\tau_{\mathbf{R}_{k,U_{m^*}},j}^u$, define the sequence of random variables to be
\begin{align}
\mathbf{Z}_{U_{m^*}}(i)=\mathbf{R}_{k,U_{m^*}}- \frac{n^{k-m^*}}{(k-m^*)!}p^{{k\choose 2}-
{m^*\choose 2}}-(\sigma-1) n^{\beta_{m^*}}\log^{\gamma_{m^*}} n.
\end{align}
In order to prove the upper bound of $\mathbf{R}_{k,U_{m^*}}$
 is the equation  in (3.5), we prove the following two claims.
\vskip 0.3cm

\noindent{\bf Claim 4.2:}\ Removing the edges of one $K_k$ in $G(i)$, we have
\begin{align*}
\mathbf{R}_{k,U_{m^*}}(i)-\mathbf{R}_{k,U_{m^*}}(i+1)=
O\bigl(n^{k-m^*-1}p^{\binom{k}{2}-\binom{m^*+1}{2}}\bigr).
\end{align*}

\begin{proof}[Proof of Claim 4.2] When we remove the edges of one $K_k$ from $G(i)$,
note that $\mathbf{R}_{k,U_{m^*}}(i)$ is the number of  $K_{k-m^*}$s in which every vertex is
in $N_{U_{m^*}}$, then  it is clearly true
that $\mathbf{R}_{k,U_{m^*}}(i)-\mathbf{R}_{k,U_{m^*}}(i+1)\geqslant 0$. Suppose the removed $K_k$ contains
 one vertex in $U_{m^*}$, denoted by $u\in U_{m^*}$;  and also contains some vertex,
 denoted by $w$, that is in $N_{U_{m^*}}$. Then the number of
 $K_{k-m^*-1}$s in which every vertex is in $N_{U_{m^*}\cup \{w\}}$
is at most $\mathbf{R}_{k,U_{m^*}\cup \{w\}}(i)$. By  the equation in (3.5),
we complete the proof.
\end{proof}

\noindent{\bf Claim 4.3:}\  The sequence $-\mathbf{Z}_{U_{m^*}}(j),-\mathbf{Z}_{U_{m^*}}(j+1),\cdots,
-\mathbf{Z}_{U_{m^*}}(\tau_{\mathbf{R}_{k,U_{m^*}},j}^u)$ is an $(\eta,N)$-bounded submartingale,
where $\eta=\Theta\bigl(n^{k-m^*-2}p^{\binom{k}{2}-\binom{m^*}{2}-1}\bigr)$
and $N=\Theta\bigl(n^{k-m^*-1}p^{\binom{k}{2}-\binom{m^*+1}{2}}\bigr)$
for $2\leqslant m^*\leqslant k-1$.

\begin{proof}[Proof of Claim 4.3]\ For all $i$ with $j\leqslant i\leqslant\tau_{\mathbf{R}_{k,U_{m^*}},j}^u$,
as the equation in (4.18), we have
\begin{equation*}
\begin{aligned}
\mathbb{E}[\Delta \mathbf{Z}_{U_{m^*}}|\mathcal{F}_i]&=\mathbb{E}[\Delta \mathbf{R}_{k,U_{m^*}}|\mathcal{F}_i]-
\frac{n^{k-m^*}}{(k-m^*)!}\Bigl[p^{{k\choose 2}-{m^*\choose 2}}(i+1)-p^{{k\choose 2}-{m^*\choose 2}}(i)\Bigr]\\
&\quad- n^{\beta_{m^*}}\log^{\gamma_{m^*}} n\bigl[\sigma(i+1)-\sigma(i)\bigr].
\end{aligned}
\end{equation*}
Note that $p=p(i)=1-k(k-1)t$, $p(i+1)=1-k(k-1)\bigl(t+ \frac{1}{n^2}\bigr)$  in (3.1),
$\sigma(i)=1- \frac{k(k-1)}{4}\log p(i)$, and $\sigma(i+1)=1-\frac{k(k-1)}{4}\log p(i+1)$
 in (3.8), then
\begin{align}
\mathbb{E}[\Delta \mathbf{Z}_{U_{m^*}}|\mathcal{F}_i]&=
\mathbb{E}[\Delta \mathbf{R}_{k,U_{m^*}}|\mathcal{F}_i]-
\frac{n^{k-m^*}}{(k-m^*)!}\biggl[-\biggl({k\choose 2}-{m^*\choose 2}\biggr)
\frac{k(k-1)}{n^2}p^{{k\choose 2}-{m^*\choose 2}-1}\notag\\
&\quad\quad+O\Bigl( \frac{1}{n^4}p^{{k\choose 2}-{m^*\choose 2}-2}\Bigr)\biggr]
-n^{\beta_{m^*}}\log^{\gamma_{m^*}} n\Bigl[ \frac{\sigma'}{n^2}+O\Bigl( \frac{\sigma''}{n^4}\Bigr)\Bigr]\notag\\
&=\mathbb{E}[\Delta \mathbf{R}_{k,U_{m^*}}|\mathcal{F}_i]+\biggl[{k\choose 2}-{m^*\choose 2}\biggr]
\frac{k(k-1)n^{k-m^*-2}}{(k-m^*)!}p^{{k\choose 2}-{m^*\choose 2}-1}\notag\\
&\quad\quad-\sigma'n^{\beta_{m^*}-2}\log^{\gamma_{m^*}} n+
O\bigl(n^{k-m^*-4}p^{{k\choose 2}-{m^*\choose 2}-2}\bigr),
\end{align}
where $O(\sigma''n^{\beta_{m^*}-4}\log^{\gamma_{m^*}} n)$  is absorbed
into $O(n^{k-m^*-4}p^{{k\choose 2}-{m^*\choose 2}-2})$ because 
$\sigma''=O(p^{-2})$ in (3.6), $\beta_{m^*}$ shown in (3.7),
$p\geqslant p_0$  in (3.9), and appropriate choices of the constants
 $\lambda$ and $\gamma_{m^*}$.

Combining the equations  in (4.17) and (4.19), we further have
\begin{align}
\mathbb{E}[\Delta \mathbf{Z}_{
U_{m^*}}|\mathcal{F}_i]
&< \frac{\bigl[{k\choose 2}-{m^*\choose 2}\bigr]k!\sigma
}{(k-m^*)!p^{{m^*\choose 2}}}n^{k-{m^*}- \frac{5}{2}}\log^{\gamma_2} n\notag\\
&\quad\quad - \frac{\bigl[{k\choose 2}-{m^*\choose 2}\bigr]k(k-1)(\sigma-1)}
{p}n^{\beta_{m^*}-2}\log^{\gamma_{m^*}} n\notag\\
&\quad\quad+ \frac{\bigl[{k\choose 2}-{m^*\choose 2}\bigr]k!\sigma(\sigma-1)
}{p^{{k\choose 2}}}n^{\beta_{m^*}-\frac{5}{2}}\log^{\gamma_2+\gamma_{m^*}} n\notag\\
&\quad\quad-\sigma'n^{\beta_{m^*}- 2}
\log^{\gamma_{m^*}} n+O\bigl(n^{k-m^*-3}p^{{k\choose 2}-{m^*\choose 2}-3}\bigr),
\end{align}
where $O(n^{k-m^*-4}p^{{k\choose 2}-{m^*\choose 2}-2})$ in (4.19)
is absorbed into $O(n^{k-m^*-3}p^{{k\choose 2}-{m^*\choose 2}-3})$ in (4.17).
At last, we  have $\mathbb{E}[\Delta \mathbf{Z}_{U_{m^*}}|\mathcal{F}_i]<0$
in (4.20)  because the following inequalities
\begin{equation*}
\begin{aligned}
& \frac{k!\sigma}{(k-m^*)!p^{{m^*\choose 2}}}n^{k-{m^*}- \frac{5}{2}}\log^{\gamma_2} n
< \frac{k(k-1)(\sigma-1)}
{2p}n^{\beta_{m^*}-2}\log^{\gamma_{m^*}} n,\\
& \frac{k!\sigma}{p^{{k\choose 2}}}n^{\beta_{m^*}- \frac{5}{2}}\log^{\gamma_2} n
< \frac{k(k-1)}
{2p}n^{\beta_{m^*}-2},\\
& O\bigl(n^{k-m^*-3}p^{{k\choose 2}-{m^*\choose 2}-3}\bigr)<
\sigma' n^{\beta_{m^*}- 2}\log^{\gamma_{m^*}} n
\end{aligned}
\end{equation*}
are obviously  true when $\beta_{m^*}$ is
in (3.7), $\sigma'= k^2(k-1)^2/4p$ in (3.8),
 $p\geqslant p_0$ in (3.9),  and appropriate choices of $\lambda$ and
 $\gamma_{m^*}$. We have proved that the
sequence $-\mathbf{Z}_{U_{m^*}}(j),-\mathbf{Z}_{U_{m^*}}(j+1),\cdots,
-\mathbf{Z}_{U_{m^*}}(\tau_{\mathbf{R}_{k,U_{m^*}},j}^u)$
is a  submartingale for any $2\leqslant m^*\leqslant k-1$.

In the following, we show the sequence is $(\eta,N)$-bounded.
By the equation in (4.18) and the calculation  in (4.19), we have
\begin{equation*}
\begin{aligned}
&-\mathbf{Z}_{U_{m^*}}(i+1)+\mathbf{Z}_{U_{m^*}}(i)\\
&=\mathbf{R}_{k,U_{m^*}}(i)- \mathbf{R}_{k,U_{m^*}}(i+1)+
\frac{n^{k-m^*}}{(k-m^*)!}\Bigl[p^{{k\choose 2}-{m^*\choose 2}}(i+1)-p^{{k\choose 2}-{m^*\choose 2}}(i)\Bigr]\\
&\quad\quad+n^{\beta_{m^*}}\log^{\gamma_{m^*}} n\bigl[\sigma(i+1)-\sigma(i)\bigr]\\
&=\mathbf{R}_{k,U_{m^*}}(i)-\mathbf{R}_{k,U_{m^*}}(i+1)-\biggl[{k\choose 2}-{m^*\choose 2}\biggr]
\frac{k(k-1)n^{k-m^*-2}}{(k-m^*)!}p^{{k\choose 2}-{m^*\choose 2}-1}\\
&\quad\quad+\sigma'n^{\beta_{m^*}-2}
\log^{\gamma_{m^*}} n+O\bigl(n^{k-m^*-4}p^{{k\choose 2}-{m^*\choose 2}-2}\bigr).
\end{aligned}
\end{equation*}
Note that $p\geqslant p_0$ in (3.9), and appropriate choices of
$\lambda$ and $\gamma_{m^*}$, we have
\begin{align*}
\biggl[{k\choose 2}-{m^*\choose 2}\biggr]
\frac{k(k-1)n^{k-m^*-2}}{(k-m^*)!}p^{{k\choose 2}-{m^*\choose 2}-1}
>\sigma'n^{\beta_{m^*}-2}
\log^{\gamma_{m^*}} n+O\bigl(n^{k-m^*-4}p^{{k\choose 2}-{m^*\choose 2}-2}\bigr).
\end{align*}
Thus, we take
\begin{align*}
\eta&=\biggl[{k\choose 2}-{m^*\choose 2}\biggr]
\frac{k(k-1)n^{k-m^*-2}}{(k-m^*)!}p^{{k\choose 2}-{m^*\choose 2}-1}\\
&=\Theta\bigl(n^{k-m^*-2}p^{{k\choose 2}-{m^*\choose 2}-1}\bigr).
\end{align*}
Since $-\mathbf{Z}_{U_{m^*}}(i+1)+\mathbf{Z}_{U_{m^*}}(i)\leqslant
\mathbf{R}_{k,U_{m^*}}(i)-\mathbf{R}_{k,U_{m^*}}(i+1)$,
applying Claim 4.2, 
we take
\begin{align*}
N=\Theta\Bigl(n^{k-m^*-1}p^{\binom{k}{2}-\binom{m^*+1}{2}}\Bigr).
\end{align*} We complete the proof of Claim 4.3.
\end{proof}

The number of the sequence $-\mathbf{Z}_{U_{m^*}}(j),-\mathbf{Z}_{U_{m^*}}(j+1),\cdots,
-\mathbf{Z}_{U_{m^*}}(\tau_{\mathbf{R}_{k,U_{m^*}},j}^u)$ is also $O(n^2p)$,
which implies $\ell=O(n^2p)$ in Lemma 2.3. Choose $a=n^{\beta_{m^*}}\log^{\gamma_{m^*}}n$, then
$a=o(\eta \ell)$. Lemma 2.3 yields that,
\begin{equation*}
\begin{aligned}
&\mathbb{P}\Bigl[\mathbf{R}_{k,U_{m^*}}\geqslant \frac{n^{k-m^*}}{(k-m^*)!}
p^{{k\choose 2}-{m^*\choose 2}}+\sigma n^{\beta_{m^*}}\log^{\gamma_{m^*}} n\Bigr]\\
&=\mathbb{P}\Bigl[-\mathbf{Z}_{U_{m^*}}(i)\leqslant -n^{\beta_{m^*}}\log^{\gamma_{m^*}} n\Bigr]\\
&\leqslant\exp\biggl[-\Omega\biggl( \frac{n^{2\beta_{m^*}}\log^{2\gamma_{m^*}} n}{n^{k-m^*}\cdot n^{k-m^*-1}}\biggr)\biggr]\\
&=\exp\biggl[-\Omega\biggl(n^{\frac{\binom{m^*}{2}-1}{\binom{k}{2}-1}}\log^{2\gamma_{m^*}} n\biggr)\biggr].
\end{aligned}
\end{equation*}
By the union bound, note that the choice to choose $j$, $m^*$ ($2\leqslant m^*\leqslant k-1$)
and $U_{m^*}\in {[n]\choose m^*}$ is at most $(k-2)n^{m^*+2}$, then
we also have
\begin{align*}
(k-2)n^{m^*+2}\exp\biggl[-\Omega\biggl(n^{\frac{\binom{m^*}{2}-1}{\binom{k}{2}-1}}\log^{2\gamma_{m^*}} n\biggr)\biggr]=o(1)
\end{align*}
because it is clearly true when $3\leqslant m^*\leqslant m-1$, and taking $\gamma_2> \frac{1}{2}$ for
$m^*=2$. In a conclusion, w.h.p.,
none of $\mathbf{R}_{k,U_{m}}$ for any $U_m\in {n\choose m}$ with $2\leqslant m\leqslant k-1$ have such a large
 upward deviations.

\begin{remark}
The argument for the lower bound of $\mathbf{R}_{k,U_{m}}$ in (3.5) for any $U_m\in {[n]\choose m}$
with $2\leqslant m\leqslant k-1$  is the symmetric analogue of
the above analysis.
\end{remark}


\section{Conclusions}

For the random $K_k$-removal algorithm, there are less direct results when $k\geqslant 4$
because their evolutionary structures are more complicated than the case $k=3$
 to investigate.  We establish dynamic concentrations
 of complete higher codegree around the expected trajectories that are derived by
 their pseudorandom properties. The final size of the random $K_k$-removal algorithm
is at most $n^{2-1/(k(k-1)-2)+o(1)}$ for $k\geqslant 4$.
In order to improve the result, it is observed that the main
 obstacle is the parameter $\mathbf{R}_{k,U_m}$ for $2\leqslant m\leqslant k-1$.
The control over
$\mathbf{R}_{k,U_m}$ loses when $p$ around $p_0$ shown in Remark 3.3,
 while the probabilities
of these extreme events are very low shown in Remark 3.4. The behaviors of these chosen random
variables $\mathbf{R}_{k,U_m}$ for $2\leqslant m\leqslant k-1$ are not in a position to analyze the structures of the process
further, and it is definitely possible to find some new ideas
 to track the random $K_k$-removal algorithm. 
 This will be investigated in future work.

\section*{Acknowledgement}\label{s:11}

Fang Tian thanks  X.-F. Pan  for helping us to point out the faults in
 the equations (3.6) and (3.7),
and some useful discussions in Remark 3.3 and 3.4. Fang Tian
was supported by the National Natural Science Foundation of China
(Grant No.~12071274). X.-F. Pan was supported by University Natural
Science Research Project of Anhui Province under Grant No. KJ2020A0001.

\section*{Appendix}

\noindent{\bf Appendix: Lower bound of $\mathbf{Q}_k(i)$ (for Remark 4.1) }
\vskip 0.3cm

\noindent For the lower bound of $\mathbf{Q}_k(i)$, we work with the critical interval
\begin{equation*}
\begin{aligned}
I_{\mathbf{Q}_k}^\ell=\Bigl( \frac{n^k}{k!}p^{{k\choose 2}}-\sigma^2n^{\alpha}p^{-1}\log^\mu n, \frac{n^k}{k!}p^{{k\choose 2}}-\sigma(\sigma-1)n^{\alpha}p^{-1}\log^\mu n\Bigr),
\end{aligned}
\eqno(1)
\end{equation*}where $\alpha$ is shown in (3.6). Consider a fixed step $j\leqslant i_0$.
Similarly, suppose $\mathbf{Q}_k(j)\in I_{\mathbf{Q}_k}^\ell$ and define
\begin{equation*}
\begin{aligned}
\tau_{\mathbf{Q}_k,j}^\ell=\min\bigl\{i_0,\max\{j,\tau\},{\mbox{the\ smallest\ }}i\geqslant j {\mbox{ such\  that}\ }\mathbf{Q}_k(i)\notin I_{\mathbf{Q}_k}^\ell\bigr\}.
\end{aligned}
\eqno(2)
\end{equation*}Let $j\leqslant i\leqslant\tau_{\mathbf{Q}_k,j}^\ell$. All calculations in this subsection
are conditioned on the estimates in (3.5) hold on
$\mathbf{R}_{k,U_m}$ for any $U_m\in {[n]\choose m}$ with $2\leqslant m\leqslant k-1$.

By the equations shown in (2.5), we get the estimate on
$\mathbb{E}[\Delta \mathbf{Q}_k| \mathcal{F}_i]$ in reverse direction,
\begin{equation*}
\begin{aligned}
\mathbb{E}[\Delta \mathbf{Q}_k| \mathcal{F}_i]
&=(-1)^{k+1}(k-1)- \frac{1}{\mathbf{Q}_k(i)}\sum_{U_2\in \mathcal{K}_2(i)}\mathbf{R}_{k,U_2}^2+\cdots+
\frac{(-1)^k(k-2)}{\mathbf{Q}_k(i)}\sum_{U_{k-1}\in \mathcal{K}_{k-1}(i)}\mathbf{R}_{k,U_{k-1}}^2\\
&>(-1)^{k+1}(k-1)- \frac{1}{\mathbf{Q}_k(i)}\Bigl( \frac{2!{k\choose 2}^2\mathbf{Q}_k^2(i)}{n^2p}+2\sigma^2 n^{2k-3}p\log^{2\gamma_2} n\Bigr)\\
&\quad\quad+  \frac{12{k\choose 3}^2\mathbf{Q}_k(i)}{n^3p}+O(n^{k-4}p^{{k\choose 2}-11})\\
&=(-1)^{k+1}(k-1)- \frac{2{k\choose 2}^2\mathbf{Q}_k(i)}{n^2p}+ \frac{12{k\choose 3}^2\mathbf{Q}_k(i)}{n^3p}+
O( n^{k-3}\sigma^2p^{-{k\choose 2}+1}\log^{2\gamma_2} n),
\end{aligned}
\end{equation*}
where $\sum_{U_2\in \mathcal{K}_2(i)}\mathbf{R}_{k,U_2}^2$ and $\sum_{U_3\in \mathcal{K}_3(i)}\mathbf{R}_{k,U_3}^2$
are replaced by the equations in (4.1) and (4.2),
the term $O(n^{k-4}p^{{k\choose 2}-11})$ comes from $\sum_{U_4\in \mathcal{K}_4(i)}\mathbf{R}_{k,U_4}^2$ in
(4.3) that dominates all the remaining terms.

Since $\mathbf{Q}_k(j)\in I_{\mathbf{Q}_k}^\ell$ shown in (1), we further have 
\begin{equation*}
\begin{aligned}
\mathbb{E}[\Delta \mathbf{Q}_k| \mathcal{F}_i]
&>(-1)^{k+1}(k-1)- \frac{2{k\choose 2}^2\bigl( \frac{n^k}{k!}p^{{k\choose 2}}-\sigma(\sigma-1)n^{\alpha}p^{-1}\log^\mu n\bigr)}{n^2p}\\
&\quad\quad+ \frac{12{k\choose 3}^2\bigl( \frac{n^k}{k!}p^{{k\choose 2}}-\sigma^2n^{\alpha}p^{-1}\log^\mu n\bigr)}{n^3p}+O\bigl( n^{k-3}\sigma^2 p^{-{k\choose 2}+1}\log^{2\gamma_2} n \bigr)\\
&=(-1)^{k+1}(k-1)- \frac{2{k\choose 2}^2 n^{k-2}}{k!}p^{{k\choose 2}-1}+ \frac{2{k\choose 2}^2\sigma(\sigma-1)n^{\alpha-2}\log^\mu n}{p^{2}}\\
&\quad\quad+ \frac{12{k\choose 3}^2n^{k-3}}{k!}p^{{k\choose 2}-1}+O\bigl( n^{k-3}\sigma^2 p^{-{k\choose 2}+1}\log^{2\gamma_2} n\bigr),
\end{aligned}
\eqno(3)
\end{equation*}where $\alpha$ is in (3.6), and
 $O( \sigma^2n^{\alpha-3}p^{-2}\log^\mu n)$
is absorbed into $O(n^{k-3}\sigma^2 p^{-{k\choose 2}+1}\log^{2\gamma_2} n)$. 

For all $i$ with $j\leqslant i\leqslant\tau_{\mathbf{Q}_k,j}^\ell$, define the sequence of random variables  to be
\begin{equation*}
\begin{aligned}
\mathbf{L}(i)=\mathbf{Q}_k(i)- \frac{n^k}{k!}p^{{k\choose 2}}+\sigma^2n^{\alpha}p^{-1}\log^\mu n.
\end{aligned}
\eqno(4)
\end{equation*}

\noindent{\bf Claim A:}\ The
sequence $\mathbf{L}(j),\mathbf{L}(j+1),\cdots,\mathbf{L}(\tau_{\mathbf{Q}_k,j}^\ell)$ is a submartingale and
the maximum one step $\Delta \mathbf{L}$ is
$O(\sigma n^{k- 5/2}\log^{\gamma_2}n)$.

\vskip 0.2cm
\begin{proof}[Proof of Claim A]\ Similarly,  for all $i$ with $j\leqslant i<\tau_{\mathbf{Q}_k,j}^\ell$,
as the equation shown in (4), we have
\begin{equation*}
\begin{aligned}
\mathbb{E}[\Delta \mathbf{L}|\mathcal{F}_i]
&=\mathbb{E}[\Delta \mathbf{Q}_k|\mathcal{F}_i]- \frac{n^k}{k!}\Bigl[p^{k\choose 2}(i+1)-p^{k\choose 2}(i)\Bigr]
+n^{\alpha}\log^\mu n\Bigl[ \frac{\sigma^2(i+1)}{p(i+1)}-\frac{\sigma^2(i)}{p(i)}\Bigr].\\
\end{aligned}
\end{equation*}
Note that $p(i)=1-k(k-1)t$,  $p(i+1)=1-k(k-1)(t+ \frac{1}{n^2})$ in (3.1),
$\sigma(i)=1- \frac{k(k-1)}{4}\log p(i)$, $\sigma(i+1)=1- \frac{k(k-1)}{4}\log p(i+1)$
 in (3.8), then by Taylor's expansion, we have
\begin{equation*}
\begin{aligned}
\mathbb{E}[\Delta \mathbf{L}|\mathcal{F}_i]
&=\mathbb{E}[\Delta \mathbf{Q}_k|\mathcal{F}_i]- \frac{n^k}{k!}\biggl[-{k\choose 2} \frac{k(k-1)}{n^2}p^{{k\choose 2}-1}
+O\Bigl( \frac{1}{n^4}p^{{k\choose 2}-2}\Bigr)\biggr]\\
&\quad\quad+n^{\alpha}\log^\mu n\biggl[ \frac{2\sigma\sigma' p-\sigma^2p'}{n^2 p^2}
+O\Bigl( \frac{\sigma^2}{n^4p^3}\Bigr)\biggr]\\
&=\mathbb{E}[\Delta \mathbf{Q}_k|\mathcal{F}_i]+ \frac{2{k\choose 2}^2 n^{k-2}
}{k!}p^{{k\choose 2}-1}+ \frac{2\sigma\sigma'n^{\alpha-2}\log^\mu n}{p}\\
&\quad\quad+ \frac{k(k-1)\sigma^2n^{\alpha-2}\log^\mu n}{p^2}+O\bigl(n^{k-4}
p^{{k\choose 2}-2}\bigr),
\end{aligned}
\eqno(5)
\end{equation*}
where  $O( n^{\alpha-4}\sigma^2 p^{-3}\log^{\mu} n)$ is absorbed
into $O(n^{k-4}p^{{k\choose 2}-2})$ because 
$\alpha$ is  in (3.6), $p\geqslant p_0$ in (3.9),
and appropriate choices of $\lambda$ and $\mu$.
Combining the equations in (3) and (5), we have
\begin{equation*}
\begin{aligned}
\mathbb{E}[\Delta \mathbf{L}|\mathcal{F}_i]
&>(-1)^{k+1}(k-1)+ \frac{\bigl[2{k\choose 2}^2+k(k-1)\bigr]\sigma^2n^{\alpha-2}\log^\mu n}{p^{2}}-
\frac{2{k\choose 2}^2\sigma n^{\alpha-2}\log^\mu n}{p^{2}}\\
&\quad\quad+ \frac{12{k\choose 3}^2n^{k-3}}{k!}p^{{k\choose 2}-1}+ \frac{2\sigma\sigma'n^{\alpha-2}\log^\mu n}{p}+
O\bigl(n^{k-3}\sigma^2 p^{-{k\choose 2}+1}\log^{2\gamma_2} n\bigr),
\end{aligned}
\end{equation*}
where $O(n^{k-4}p^{{k\choose 2}-2})$ in (5) is absorbed into
$O(n^{k-3}\sigma^2 p^{-{k\choose 2}+1}\log^{2\gamma_2} n)$ in (3).
We have
\begin{align*}
2\sigma\sigma'n^{\alpha-2}\log^\mu np^{-1}= 2{k\choose 2}^2\sigma n^{\alpha-2}p^{-2}\log^\mu n
\end{align*}
by  $\sigma'= \frac{1}{4}k^2(k-1)^2 p^{-1}$ in (3.8).
It follows that
\begin{equation*}
\begin{aligned}
\mathbb{E}[\Delta \mathbf{L}|\mathcal{F}_i]
&>(-1)^{k+1}(k-1)+ \frac{\bigl[2{k\choose 2}^2+k(k-1)\bigr]\sigma^2n^{\alpha-2}\log^\mu n}{p^{2}}\\
&\quad\quad+ \frac{12{k\choose 3}^2n^{k-3}}{k!}p^{{k\choose 2}-1}+
O\bigl(n^{k-3}\sigma^2 p^{-{k\choose 2}+1}\log^{2\gamma_2} n\bigr).
\end{aligned}
\eqno(6)
\end{equation*}
Note that
\begin{align*}
\Bigl[2{k\choose 2}^2+k(k-1)\Bigr]\sigma^2n^{\alpha-2}p^{-2}\log^\mu n>
O(n^{k-3}\sigma^2 p^{-{k\choose 2}+1}\log^{2\gamma_2} n)+(-1)^{k+1}(k-1)
\end{align*}
when $\alpha$ is in (3.6), $p\geqslant p_0$ in (3.9),
and appropriate choices of $\lambda$, $\mu$ and $\gamma_2$. 
We have $\mathbb{E}[\Delta \mathbf{L}|\mathcal{F}_i]>0$.
The sequence $\mathbf{L}(j),\mathbf{L}(j+1),\cdots,\mathbf{L}(\tau_{\mathbf{Q}_k,j}^\ell)$ is a submartingale.

Next, we show the maximum one step $\Delta \mathbf{L}$ is
$O(\sigma n^{k- 5/2}\log^{\gamma_2} n)$. As the equation in (4) and the calculation in (5),
we have
\begin{equation*}
\begin{aligned}
\Delta \mathbf{L}&=\Delta \mathbf{Q}_k+ \frac{2{k\choose 2}^2 n^{k-2}
}{k!}p^{{k\choose 2}-1}+ \frac{2\sigma\sigma'n^{\alpha-2}\log^\mu n}{p}+
\frac{k(k-1)\sigma^2n^{\alpha-2}\log^\mu n}{p^2}+O(n^{k-4}p^{{k\choose 2}-2}).
\end{aligned}
\end{equation*}Apply the equation of $\Delta \mathbf{Q}_k$ in (2.4)
and the estimates on $\mathbf{R}_{k,U_m}$ for any $U_m\in {[n]\choose m}$
when $2\leqslant m\leqslant k-1$ in (3.5) to the above display, then
\begin{equation*}
\begin{aligned}
\Delta \mathbf{L}&\leqslant -{k\choose 2}\Bigl( \frac{n^{k-2}}{(k-2)!}p^{{k\choose 2}-1}
-\sigma n^{\beta_2}\log^{\gamma_2} n\Bigr)
+{k\choose 3}\Bigl( \frac{n^{k-3}}{(k-3)!}p^{{k\choose 2}-{3\choose 2}}+
\sigma n^{\beta_3}\log^{\gamma_3} n\Bigr)+\cdots\\
&\quad\quad+ \frac{2{k\choose 2}^2 n^{k-2}}{k!}p^{{k\choose 2}-1}+
\frac{2\sigma\sigma'n^{\alpha-2}\log^\mu n}{p}+ \frac{k(k-1)\sigma^2n^{\alpha-2}\log^\mu n}{p^2}
+O\bigl(n^{k-4}p^{{k\choose 2}-2}\bigr)\\
&=O(\sigma n^{k- \frac{5}{2}}\log^{\gamma_2} n),
\end{aligned}
\end{equation*}where $\alpha$ and
$\beta_2$ are in (3.6) and (3.7), the term
$O(\sigma^2n^{\alpha-2}p^{-2}\log^\mu n)$ is absorbed into $O(\sigma n^{k- 5/2}\log^{\gamma_2} n)$
when $p\geqslant p_0$ shown in (3.9), and
appropriate choices of $\lambda$, $\mu$
and $\gamma_2$.

The number of steps in this sequence is also 
 $O(n^2p)$.
Since $\mathbf{Q}_k(j)\in I_{\mathbf{Q}_k}^\ell$
 in (1), we have $\mathbf{L}(j)<\sigma n^{\alpha}p^{-1}\log^\mu n$
from (4). For all $i$ with
$j\leqslant i\leqslant\tau_{\mathbf{Q}_k,j}^\ell$, Lemma 2.3 yields that the probability of such a large
deviation beginning at the step $j$ is at most
\begin{equation*}
\begin{aligned}
&\mathbb{P}\Bigl[\mathbf{Q}_k(i)\leqslant
\frac{n^k}{k!}p^{{k\choose 2}}-\sigma^2n^{\alpha}p^{-1}\log^\mu n\Bigr]\\
&=\mathbb{P}\Bigl[\mathbf{L}(i)\leqslant 0\Bigr]\\
&\leqslant \exp\biggl[-\Omega\biggl( \frac{\bigl(\sigma n^{\alpha}p^{-1}\log^\mu n\bigr)^2}
{(n^2p)\bigl( \sigma n^{k- 5/2}\log^{\gamma_2}n\bigr)^2}\biggr)\biggr]\\
&=\exp\biggl[-\Omega\biggl(\frac{n^{2\alpha-2k+3}\log^{2\mu} n}{p^3\log^{2\gamma_2} n}\biggr)\biggr].
\end{aligned}
\end{equation*}
By the union bound, note that there are at most $n^2$ possible values of $j$
 shown in (3.1), then we have
 \begin{equation*}
n^2\exp\biggl[-\Omega\biggl(\frac{n^{2-\frac{2}{\binom{k}{2}-1}}\log^{2\mu} n}{p^3\log^{2\gamma_2} n}\biggr)\biggr]=o(1)
\end{equation*}
 with $\alpha$ is shown in (3.6). W.h.p. $\mathbf{Q}_k(i)$ never crosses
its critical interval $I_{\mathbf{Q}_k}^\ell$  in (4),
 and so the lower bound on $\mathbf{Q}_k(i)$ in (3.4) is true.
\end{proof}

\end{document}